\documentclass{amsart}

\usepackage[normalem]{ulem}

\usepackage{lineno}

\usepackage{caption}

\usepackage{soul}

\usepackage{graphicx}

\usepackage{amssymb, amsmath, amsthm, hyperref, euscript, color}
\usepackage[dvipsnames]{xcolor}

\usepackage{amscd}

\usepackage{bbm}

\usepackage{enumitem}

\usepackage[english]{babel}

\numberwithin{equation}{section}

\theoremstyle{plain}

\newtheorem{theorem}{Theorem}

\newtheorem{corollary}[theorem]{Corollary}

\newtheorem{remark}[theorem]{Remark}

\newtheorem{question}{Question}

\newtheorem{thm}{Theorem}

\theoremstyle{definition}

\newtheorem{definition}[theorem]{Definition}

\newcommand{\Z}{\mathbb{Z}}

\newcommand{\N}{\mathbb{N}}

\newcommand{\C}{\mathbb{C}}

\DeclareMathOperator{\Kod}{Kod}

\makeatletter
\@namedef{subjclassname@2020}{%
  \textup{2020} Mathematics Subject Classification}
\makeatother

\begin{document}

\author{Rafael Torres}

\title[Stable generalized complex structures on sphere bundles]{An observation on the existence of stable generalized complex structures on ruled surfaces.}


\address{Scuola Internazionale Superiori di Studi Avanzati (SISSA)\\ Via Bonomea 265\\34136\\Trieste\\Italy}

\email{rtorres@sissa.it}


\maketitle

\emph{Abstract}:  We point out that any stable generalized complex structure on a sphere bundle over a closed surface of genus at least two must be of constant type. 

\section{Introduction.}\label{Introduction}

The canonical line bundle $K_{\mathcal{J}}\subset \wedge^{\bullet}T^\ast_{\C}M$ of a stable generalized complex structure $\mathcal{J}$ on a smooth n-manifold $M$ is generated pointwise by the complex differential form\begin{equation}\label{Local}\rho = e^{B + i\omega}\wedge \Omega,\end{equation}where $B$ and $\omega$ are real two-forms and the complex form on the right-side of (\ref{Local}) decomposes into $\Omega = \theta_1\wedge\cdots \wedge \theta_k$ for $(\theta_1, \ldots, \theta_k)$ a basis for $L\cap T^\ast_{\C}$ and $L$ is the $+i$-eigenbundle of $\mathcal{J}$  (\cite{[GualtieriPhD], [Gualtieri]} and Cavalcanti-Gualtieri \cite{[CavalcantiGualtieri3]} for background and further details). A fundamental invariant of a stable generalized complex structure $(M, \mathcal{J})$ is its type and its defined by\begin{center}Type($\mathcal{J})$ = $\deg(\Omega) = k$;\end{center}see \cite[Section 3.1]{[Gualtieri]}. In dimension four, the possible constant types for a given $(M, \mathcal{J})$ are $k\in \{0,1, 2\}$. The value $k = 0$ corresponds to a stable generalized complex structure induced by a symplectic structure, while $k = 2$ corresponds to one induced by a complex structure. Stable generalized complex structures of type $k = 1$ have been studied by Bailey-Cavalcanti-Gualtieri \cite{[BaileyCavalcantiGualtieri]} and Chen-Nie \cite{[ChenNie]}.

Moreover, there is now a myriad of examples of stable generalized complex structures in dimension four whose type jumps from $k = 0$ to $k = 2$ \cite{[Gualtieri], [CavalcantiGualtieri1], [CavalcantiGualtieri2], [CavalcantiKlaasseWitte], [GotoHayano], [Torres], [TorresYazinski]}. As Cavalcanti-Gualtieri have shown, for these kind of stable generalized complex structures the underlying smooth four-manifold need not support a symplectic nor a complex structure \cite{[CavalcantiGualtieri1], [CavalcantiGualtieri2]}.  The type change occurs along a smoothly embedded torus $T\subset M$, which is a path connected-component of the type change locus of $\mathcal{J}$. All known examples of such stable generalized complex structures in the literature share the common trait that the underlying closed four-manifold has nonnegative Euler characteristic. This raises the following question.

\begin{question}\label{Question A}Does there exist a closed four-manifold of negative Euler characteristic that admits a stable generalized complex structure with non-empty type change locus?
\end{question}

The main result of this short note puts Question \ref{Question A} into perspective by exhibiting four-manifolds of negative Euler characteristic that solely admit stable generalized complex structures of constant type, a behaviour which was not previously known to occur; cf. Remark \ref{Remark Examples}.

\begin{thm}\label{Theorem A} The type change locus of a stable generalized complex structure on a sphere bundle over a closed surface of genus strictly greater than one is empty. 

In particular, any stable generalized complex structure on such a four-manifold arises from a K\"ahler structure. 
\end{thm}

This note aims to make well-known classification theorems and results addressing the minimal genus function of a symplectic four-manifold (see \cite[Chapter 2]{[GompfStipsicz]}) more widespread within the generalized complex geometry community. The main ingredient in the proof of Theorem \ref{Theorem A} is Li-Li's solution to the minimal genus problem of smooth embeddings in sphere bundles over surfaces in \cite{[BLiLi]} (see Theorem \ref{Theorem LL1} and Theorem \ref{Theorem LL2}). The other ingredients in the proof are a pair of generalized complex cut-and-paste constructions of Cavalcanti-Gualtieri \cite{[CavalcantiGualtieri1], [CavalcantiGualtieri2]}, the symplectic Thom conjecture due to Oszv\'ath-Szab\'o \cite[Theorem 1.1]{[OS]}, and Liu's classification of irrational ruled symplectic four-manifolds \cite{[Liu]} (see Theorem \ref{Theorem Diffeo Type}).

\subsection{Acknowledgements:}This note arose out of a conversation with Gil Cavalcanti during the BIRS workshop "Generalized geometry meets string theory" at the Instituto de Matem\'aticas de la Universidad de Granada (IMAG). I thank the organizers for their hospitality. I thank the referee for several useful suggestions to improve the manuscript.

\section{Background results.}\label{Background results}

Key background results that are used in the proof of Theorem \ref{Theorem A} are collected in this section with the purpose of making the short note as self-contained as possible.

\subsection{A relation between the symplectic Kodaira dimension and a stable generalized complex surgery}Cavalcanti-Gualtieri introduced in \cite[Section 3]{[CavalcantiGualtieri1]} and \cite[Section 4]{[CavalcantiGualtieri2]} a four-dimensional cut-and-paste operation that has been a fruitful source of examples of stable generalized complex four-manifolds; cf. Goto-Hayano \cite{[GotoHayano]}. Their main results can be rephrased into the following theorem. 

\begin{theorem}\label{Theorem CG} Cavalcanti-Gualtieri \cite{[CavalcantiGualtieri1], [CavalcantiGualtieri2]}.Let $(\widehat{M}, \mathcal{J})$ be a stable generalized complex four-manifold whose type change locus has path-connected components $\{\widehat{T}_1, \ldots, \widehat{T}_n\}$ for some $n\in \N$. Suppose that each of these tori has self-intersection zero. There is a symplectic four-manifold $(M, \omega)$ with Euler characteristic and signature given by\begin{center}$\chi(\widehat{M}) = \chi(M)$ and $\sigma(\widehat{M}) = \sigma(M)$,\end{center}as well as symplectically embedded tori $T_1, \ldots, T_n\subset (M, \omega)$ of self-intersection zero. 
\end{theorem}

The second Stiefel-Whitney class of the four-manifolds of Theorem \ref{Theorem CG} need not coincide. Indeed, the setting of the example considered by Cavalcanti-Gualtieri \cite[Example 4.2]{[CavalcantiGualtieri1]} is $\widehat{M} = 3\mathbb{CP}^2\# 19\overline{\mathbb{CP}^2}$ and $M$ is a K3 surface. 






Before discussing a relation between the symplectic Kodaira dimension of $M$ and the characteristic numbers of $\widehat{M}$, first recall the definition of the former. Let $(M, \omega)$ be a minimal symplectic four-manifold and let $K_\omega$ be its canonical class. The symbol $K_\omega$ also denotes the first Chern class of any almost-complex structure on $M$ that is compatible with the symplectic structure $\omega$.

\begin{definition}\label{Definition Kodaira}The symplectic Kodaira dimension. Li \cite[Definition 2.2]{[Li]}. The symplectic Kodaira dimension $\Kod(M, \omega)$ of a minimal symplectic four-manifold $(M, \omega)$ is defined by\[
\Kod(M, \omega) =
\begin{cases}
-\infty & \text{if}\enspace K_\omega\cdot [\omega] < 0 \; \text{or} \enspace K_\omega \cdot K_\omega < 0,\\
0 & \text{if}\enspace K_\omega\cdot [\omega] = 0 \enspace \text{and} \enspace K_\omega \cdot K_\omega = 0,\\
1 & \text{if} \enspace K_\omega\cdot [\omega] > 0\enspace  \text{and} \enspace K_\omega \cdot K_\omega = 0 \enspace\text{or}\\
2 & \text{if}\enspace K_\omega\cdot [\omega] > 0 \enspace \text{and} \enspace K_\omega \cdot K_\omega > 0. 
\end{cases}
\]

The symplectic Kodaira dimension of a non-minimal symplectic four-manifold is defined as the symplectic Kodaira dimension of any of its minimal models. 
\end{definition}

A minimal symplectic four-manifold $(M, \omega)$ with $\pi_1(M)$ of exponential growth has $\Kod(M, \omega) = - \infty$ if and only if $M$ is irrationally ruled by a result of Liu \cite{[Liu]}. All such four-manifolds admit a K\"ahler structure  and are  diffeomorphic to either the product $S^2\times \Sigma_g$ or the non-trivial bundle $S^2\widetilde{\times} \Sigma_g$ in the abscence of smoothly embedded spheres of self-intersection minus one; see \cite[Section 3]{[Li]} for a guide to the literature on the beautiful characterizations of these four-manifolds. For our purposes, we state Liu's classification as the following result.

\begin{theorem}\label{Theorem Diffeo Type} Liu \cite[Main Theorem A, Theorem B]{[Liu]}. Let $(M, \omega)$ be a symplectic four-manifold with $\Kod(M, \omega) = - \infty$ and whose fundamental group has exponential growth. Then, $M$ is diffeomorphic to $M_g\#k\overline{\mathbb{CP}^2}$ where $M_g$ is a sphere bundle over a surface of genus at least two and $k\in \Z_{\geq 0}$.
\end{theorem}

The Euler characteristic of a minimal symplectic four-manifold $(M, \omega)$ with $\Kod(M, \omega) = - \infty$ is negative if and only if its fundamental group has exponential growth. The precise diffeomorphism type of the four-manifold of Theorem \ref{Theorem Diffeo Type} is determined from its intersection form over $\Z$ \cite[Definition 1.2.1]{[GompfStipsicz]} and/or its Euler characteristic,  signature and second Stiefel-Whitney class. 






\subsection{The minimal genus of symplectic surfaces in irrational ruled symplectic four-manifolds}\label{Section Genus}The main ingredient in the proof of Theorem \ref{Theorem A} is work of B.-H. Li and T.-J. Li \cite[Theorems 1 and 2]{[BLiLi]} that solves the minimal genus problem for smooth embeddings of surfaces in irrational ruled symplectic four-manifolds. We now recall their statements in order to make this note as self-contained as possible. 

\begin{theorem}\label{Theorem LL1}Li-Li \cite[Theorem 1]{[BLiLi]}. Let $a$ and $b$ be nonnegative integers, and consider the second homology classes\begin{center}$x = [\{pt\}\times \Sigma_g]$ and $y = [S^2\times \{pt\}]$ .\end{center}The minimal genus $g_{\xi}$ of $\xi = a x + b y\in H_2(S^2\times \Sigma_g)$ is given by

\[
g_{\xi} =
\begin{cases}
(a - 1)(g - 1 + b) + g &   (b + g)a\neq 0\\
0 & \text{otherwise}.
\end{cases}
\]

\end{theorem}

The two orientable sphere bundles over a closed surface $\Sigma_g$ are distinguished by their second Stiefel-Whitney class. The corresponding statement in the case where the second Stiefel-Whitney class does not vanish  is as follows. 

\begin{theorem}\label{Theorem LL2}Li-Li \cite[Theorem 2]{[BLiLi]}. Let $S^2\widetilde{\times} \Sigma_g$ be the non-trivial sphere bundle over a surface of genus $g > 0$. Let $x$ be the homology class of its section and $y$ the homology class of a fiber. The minimal genus $g_{\xi}$ of $\xi = ax + by\in H_2(S^2\widetilde{\times} \Sigma_g)$ is given by\[
g_{\xi} =
\begin{cases}
(a'- 1)(g - 1 + \frac{1}{2}|a' + 2b'|) + g &   a\neq 0\\
0 & \text{otherwise}
\end{cases}
\]where $a' = |a|$ and $b' = \frac{a}{|a|}b$ when $a\neq0$.

\end{theorem}

The following consequence of Theorem \ref{Theorem LL1} and Theorem \ref{Theorem LL2} was pointed out by the referee.

\begin{corollary}\label{Corollary Referee}Let $M$ be a sphere bundle over a surface of genus strictly greater than one. If a homology class $\xi \in H_2(M)$ is represented by a torus $T$, then $\xi$ is a multiple of the homology class of the fiber of $M$. In particular, $\xi^2 = 0$ and $T$ is not a minimum genus representative of $\xi$.
\end{corollary}


\section{Proof of Theorem \ref{Theorem A}}\label{Section Proof Theorem A}

There are stable generalized complex structures of constant type $k\in \{0, 2\}$ on both bundles $S^2\times \Sigma_g$ and $S^2\widetilde{\times} \Sigma_g$ that come from their structure as irrational ruled complex surfaces. They can also be quipped with stable generalized complex structures of constant type $k = 1$ by considering the sphere as the complex projective line and the volume form on the closed surface of genus $g$.

To see that neither of these four-manifolds admits a stable generalized complex structure with at least one connected component in its type change locus, we proceed by contradiction. Let $\widehat{M}$ be a sphere bundle over a closed surface of genus at least two and suppose that there is a stable generalized complex structure $(\widehat{M}, \mathcal{J})$ whose type change locus has at least one connected component. The following two scenarios must be considered. 
\begin{enumerate}
\item There is at least one path-connected component $\widehat{T}\subset \widehat{M}$ of the type change locus of $\mathcal{J}$ with non-zero self-intersection.
\item Every path-connected component $\widehat{T}$ of the type change locus of $\mathcal{J}$ has self-intersection zero.
\end{enumerate}

The first scenario is immediately ruled out by Corollary \ref{Corollary Referee}. In order to rule out the second scenario, we apply Theorem \ref{Theorem CG} to $(\widehat{M}, \mathcal{J})$ to obtain a symplectic four-manifold $(M, \omega)$ with signature $\sigma(M) = \sigma(\widehat{M}) = 0$ and Euler characteristic $\chi(M) = \chi(\widehat{M}) < 0$. These values imply that $(M, \omega)$ is a minimal symplectic four-manifold whose fundamental group has exponential growth. Theorem \ref{Theorem CG} implies that there is a symplectic torus $T\subset M$ of self-intersection zero, where $M$ is either $S^2\times \Sigma_g$ or $S^2\widetilde{\times} \Sigma_g$ by Liu's classification of irrational ruled symplectic four-manifolds in Theorem \ref{Theorem Diffeo Type}; cf. \cite[Theorem 10.1.19]{[GompfStipsicz]}. Such symplectic surface would be genus minimizing by the proof of the symplectic Thom conjecture of Oszv\'ath-Szab\'o \cite[Theorem 1.1]{[OS]}. However, this contradicts Corollary \ref{Corollary Referee}. We conclude that the second scenario is not possible.

\hfill $\square$

\begin{remark}\label{Remark Examples} It is well-known that the number of path-connected components of the type change locus of a stable generalized complex structure on $S^2\times S^2$ and $S^2\times T^2$ can be chosen to be arbitrarily large \cite{[Torres], [TorresYazinski], [GotoHayano]}. 

\end{remark}

\end{document}